\documentclass[reqno,13pt]{amsart}
\usepackage{bbm}
\def\1{\mathbbm{1}}

\usepackage{amsmath,amsthm}
\usepackage[mathscr]{euscript}
\usepackage{graphicx}
\usepackage{amsmath,amssymb}
\usepackage{circuitikz}
\usepackage{geometry} 

\usepackage{pstricks}
\usepackage{pst-node}

\usepackage{tikz}
\usetikzlibrary{arrows,chains,matrix,positioning,scopes}
\tikzstyle{element}=[rectangle,draw,fill=white, line width=1pt]
\tikzstyle{terminal}=[circle,draw, scale=0.3, line width=1pt,red]
\tikzstyle{fleche}=[->,>=stealth', very thick]
\tikzstyle{fleche1}=[->,>=stealth', very thick, red]
\usepackage{placeins}
\definecolor{darkgoldenrod4}{rgb}{0.55,0.4,0.55}
\definecolor{maroon4}{rgb}{0.55,0.11,0.38}
\definecolor{indianred}{rgb}{0.8,0.36,0.36}
\definecolor{purple1}{rgb}{0.61,0.19,1}
\definecolor{goldenrod1}{rgb}{1,0.76,0.15}
\definecolor{indianred3}{rgb}{0.8,0.33,0.33}
\definecolor{red4}{rgb}{0.55,0,0}
\definecolor{darkslategray}{rgb}{0.18,0.31,0.31}
\definecolor{firebrick}{rgb}{0.7,0.13,0.13}
\definecolor{slateblue3}{rgb}{0.41,0.35,0.8}
\definecolor{mediumorchid4}{rgb}{0.48,0.22,0.55}
\definecolor{thistle4}{rgb}{0.55,0.48,0.55}
\definecolor{rltred}{rgb}{0.75,0,0}
\definecolor{rltgreen}{rgb}{0,0.5,0}
\definecolor{oneblue}{rgb}{0,0,0.75}
\definecolor{marron}{rgb}{0.64,0.16,0.16}
\definecolor{forestgreen}{rgb}{0.13,0.54,0.13}
\definecolor{purple}{rgb}{0.62,0.12,0.94}
\definecolor{dockerblue}{rgb}{0.11,0.56,0.98}
\definecolor{freeblue}{rgb}{0.25,0.41,0.88}
\definecolor{myblue}{rgb}{0,0.2,0.4}

\def\R{{\mathbb{R}}}

\def\la{\lambda}

\def\t{\tau}

\def\calL{{\mathcal{L}}}

\def\calS{{\mathcal{S}}}

\def\R{\mathbb R}

\def\S{\mathbb S}
\def\C{\mathbb C}

\def\M{\mathbb M}

\def\T{\mathbb T}

\newtheorem{definition}{Definition}[section]

\newtheorem{lemma}{Lemma}[section]
\newtheorem{proposition}{Proposition}[section]
\newtheorem{corollary}{Corollary}[section]
\newtheorem{theorem}{Theorem}[section]
\theoremstyle{remark}
\newtheorem{remark}{Remark}[section]
\newtheorem{example}{Example}[section]

\title{Admissibility of control operators for positive semigroups and robustness of input-to-state stability} 

\author{Yassine El Gantouh} \address{Yassine El Gantouh, School of Mathematical Sciences, Zhejiang Normal University, Jinhua, P.R. China}\email{elgantouhyassine@gmail.com}

\author{Yang Liu} \address{Yang Liu, School of Mathematical Sciences, Zhejiang Normal University, Jinhua, P.R. China}\email{liuyang@zjnu.edu.cn}

\author{Jianquan Lu} \address{Jianquan Lu, School of Mathematics, Southeast University, Nanjing, P.R. China}\email{jqluma@seu.edu.cn} 

\author{Jinde Cao} \address{Jinde Cao, Jiangsu Provincial Key Laboratory of Networked Collective Intelligence, Nanjing, P.R. China}\email{jdcao@seu.edu.cn}

\subjclass[2020]{93C05. 93C28. 47B65. 46B42. 37C75} 
\keywords{Distributed parameter systems. Admissible control operator. Positive linear operators and order-bounded operators. Input-to-state stability. Robustness}

\begin{document}
	\maketitle
	
	\renewcommand{\sectionmark}[1]{}
	\begin{abstract}
    In this paper, we investigate well-posedness and stability properties of distributed parameter systems, with particular emphasis on linear positive control systems. We establish a characterization of the well-posedness in the Banach lattice setting. Furthermore, we derive a resolvent condition for admissibility of control operators for positive semigroups. In addition, we study the behavior of input-to-state stability (ISS) under unbounded perturbations of the underlying semigroup generator. More precisely, we establish necessary and sufficient conditions for the robustness of ISS under Desch-Schappacher perturbations. Our theoretical results are demonstrated through a boundary value-controlled transport equation with non-local boundary conditions.
	\end{abstract}

\section{Introduction}\label{Sec:1}
In this paper, we investigate the admissibility of input operators and the robustness of ISS under Desch-Schappacher perturbations. To specify our terminology and notation, we consider the abstract differential equation
\begin{align}\label{diff-equ}
	\dot{z}(t) = A z(t)+B u(t),\quad t> 0,\qquad z(0)=x,
\end{align}
where $ A$ generates a C$_0$-semigroup $ \mathbb{T}:=(T(t))_{t\ge 0}$ on $ X $, a Banach space with norm $\Vert \cdot\Vert$. Let $ \rho (A) $ denote the resolvent set of $ A $ and let $ R(\la,A): =(\la I_X-A)^{-1} $ denote the resolvent operator of $ A $. The growth bound $\omega_{0}(A)$ of $\T$ is given by $ \omega_{0}(A):=\inf \lbrace t^{-1}\log \Vert T(t)\Vert: t> 0\rbrace $. The spectral bound $s(A)$ of $A$ is defined by $s(A):=\sup\{{\rm Re}\, \la : \la \in \sigma(A)\}$, where $\sigma(A):= \C\backslash \rho (A)$ is the spectrum of $A$. We let $ X_{-1,A} $ denote the completion of $X$ in the norm $\Vert x\Vert_{-1}:=\Vert R(\la ,A)x\Vert $ for some $ \la \in \rho (A) $ (see, e.g., \cite[Sect. II.5]{NEN}). The control space $U$ is a Banach space with norm $\Vert \cdot\Vert_U$, and $B:U\to X_{-1,A}$ is a control or input operator. In this context, $X$ is referred to as the state-space and $U$ as the input-space. It is well-known (see, e.g., \cite{TW}) that the mild solutions of \eqref{diff-equ} are given by 
\begin{align}\label{S2.2}
	z(t,x,u)&=T(t)x+\int_{0}^{t}T_{-1}(t-s)Bu(s)ds,
\end{align}
for all $t\ge 0$, $ x\in X $, and $u\in  L^{p}_{loc}(\R_+;U) $. The integral in \eqref{S2.2} is calculated in $ X_{-1,A} $. For practical purposes, we seek solutions $z$ that are continuous $X$-valued functions, motivating the following definition of admissibility (see, e.g., \cite[Def. 4.1]{WC}).
\begin{definition}\label{S2.D1}
	Let $A$	generate a C$_0$-semigroup $\T$ and $p\in [1,+\infty]$. An operator $ B\in \mathcal{L}(U,X_{-1,A}) $ is called a $ L^p $-admissible control operator for $A$ (or $(A,B)$ is $L^p$-admissible for short) if, for some (and hence for all) $\t>0$ and all $ u\in L^{p}(\R_+;U) $,
	\begin{align}\label{S2.2'}
		\Phi_{\t}u:=\int_{0}^{\t}T_{-1}(\t-s)Bu(s)ds\in X.
	\end{align}
\end{definition}
\vspace{.2cm}

Key consequences of admissibility include:
\begin{itemize}
	\item For every $t\ge 0$ there exists a constant $\kappa:=\kappa(t)>0 $ such that
	\begin{align}\label{input-map}
		\Vert \Phi_{t}u\Vert \leq \kappa\Vert u\Vert_{ L^{p}([0,t];U)},
	\end{align}
	for all $ u\in L^{p}(\R_+;U) $;
	\item For each $\alpha>\omega_{0}(A)$ there exits $m_\alpha>0$ such that 
	\begin{align}\label{EstimateB}
		\Vert R(\la,A_{-1})B\Vert \le \dfrac{m_\alpha}{\sqrt[p]{{\rm Re}\, \la-\alpha}},\qquad\forall\, {\rm Re}\, \la>\alpha;
	\end{align}
	\item The solutions $z(\cdot,x,u)$ are continuous $X$-valued functions for any $f\in X$ and $u\in L^p_{loc}(\R_+,U)$ where $p\in [1,+\infty)$, cf. \cite[Prop. 2.3]{WC} (see also \cite[Sect. 2.2]{Salam} for $p=2$ and \cite[Prop 2.5]{JaNPS} for $p=\infty$). 
\end{itemize}
\vspace{.2cm}

In particular, the admissibility of input operators is a basis for system analysis. This fact has motivated several works to establish sufficient conditions that guarantee the admissibility of $B$ and are easy to verify in practice. We point out here that the resolvent condition \eqref{EstimateB} is a necessary condition for $ L^p $-admissibility, and its sufficiency has been established (for $p=2$) in \cite{Wconj,Weiss} for certain class of semigroups. In these papers, the author also conjectured that $L^2$-admissibility and the resolvent condition \eqref{EstimateB} are equivalent when $X$ and $U$ are Hilbert spaces. However, this conjecture was disproven, and counterexamples were given in \cite{JaZ,Wy,ZJS}. On the other hand, in \cite{Hao}, the authors established a sufficient condition for $L^2$-admissibility of the input operator for diagonal semigroups on $\ell^2$ with scalar input. In \cite{Wl} it was shown that the condition in \cite{Hao} is also necessary for admissibility. An extension of the above results to the case $\ell^p$ ($1\le p<\infty$) was given in \cite{MU}. In \cite{Engel}, the author gave a semigroup criterion for $L^p$-admissibility, the origin of which can be traced back to \cite{Grab}. The robustness of $L^p$-admissibility under perturbations of the underlying semigroup generator was discussed in \cite{Hadd,HaWe,Mei,Schmid} (see also \cite[Sect. 5.5]{TW}). For a more comprehensive exploration of criteria for admissibility of input operators, readers are directed to the survey paper \cite{JaP} and the books \cite{Staf,TW}.
\vspace{.1cm}

In many PDEs models some constraints need to be imposed when considering concrete applications. This is for instance the case of evolutionary networks (heat conduction, transportation networks, etc.) where realistic models have to take into account that the state and/or the input must satisfy some positivity constraints to guarantee its physical relevance (see, e.g., \cite{Colombo,El,LTZ,Tilman}). 
\vspace{.1cm}

In this paper, we investigate the positivity property of \eqref{diff-equ} and admissibility of control operators for positive semigroups. Firstly, we derive a characterization for the well-posedeness of positive control systems. In particular, we show that the admissibility of a positive control operator for a given positive semigroup is completely determined by its behaviour with respect to positive inputs. Then, we establish a sufficient condition for $L^1$-admissibility in non-reflexive Banach spaces. This condition is based on a lower norm estimate for the resolvent operator on the positive cone. This resolvent condition finds its origin in problems of semigroup generation (\cite{BaRo,Arendt}), which has already found several applications \cite{BT,El1,El2,LaMo}. Secondly, our investigation delves into ISS property with respect to unbounded positive perturbations of the underlying semigroup generator. Notably, we establish necessary and sufficient conditions for the robustness of ISS under the class of Desch-Schappacher perturbations (see Section \ref{Sec:4}). Typically, these perturbations appear in many important applications, for instance, in population dynamic models, delay differential equations and generally in boundary Cauchy problems. We refer to \cite{Desh} and \cite[Sect. III.3]{NEN} for abstract results concerning this kind of perturbations and \cite{AGPF,BJVW,El1,Wintermayr} for recent applications in context of positive dynamical systems. 
\vspace{.1cm}

For related results on admissibility of (positive) input operators for positive semigroups we refer to \cite{AGPF,Wintermayr}. Notably, in \cite[Sect. 4]{AGPF}, the authors proved that admissibility can be deduced from the order and geometric properties of the state space and/or input space.
\vspace{.1cm}

Outline of the paper: In Section \ref{Sec:2}, we first discuss well-posedness and positive property of \eqref{diff-equ}. Then, we establish a sufficient condition for $L^1$-admissibility of control operators for positive semigroups (Theorem \ref{S2.T1}). Section \ref{Sec:4} is devoted to study the behavior of ISS under Desch-Schappacher perturbations. In Section \ref{Sec:5}, we illustrate our framework by an application to a boundary controlled transport equation with non-local boundary condition.
\vspace{.2cm}

{\bf Notation and terminology:} Throughout this paper $\mathbb{C},\mathbb{R},\mathbb{R}_+,\mathbb{N}$ denote the sets of complex numbers, real numbers, positive real numbers, and natural numbers, respectively. Let $E$ and $F$ be Banach spaces. We denote by $\calL(E,F)$ the space of bounded linear operators from $E$ to $F$, equipped with the operator norm, and we write $\calL(E):=\calL(E,E)$ for the algebra of bounded linear operators on $E$. For an operator $P\in \mathcal{L}(E,F)$ and a subspace $Z$ of $E$, we denote by $P\vert_{Z}$ the restriction of $P$ on $Z$, and the range of $P$ is denoted by Range$(P)$. For a real $p\ge 1$, we denote by $L^p(\R_+;E)$ the space of all measurable $p$-integrable functions $f:\R_+\to E$.
\vspace{.1cm}

A Banach lattice $E:=(E,\le)$ is a partially ordered Banach space for which any given elements $x,y$ of $E$ have a supremum $x \vee y$ and for all $x,y,z\in E$ and $\alpha\ge 0$
\begin{align*}
	\left\lbrace
	\begin{array}{lll}
		x\leq y &\implies (x+z\le y+z \quad {\rm and} \quad \alpha x\leq \alpha y ),\\
		\vert x\vert \leq \vert y\vert & \implies \Vert x\Vert \leq \Vert y\Vert,
	\end{array}
	\right.
\end{align*}
with, for all $x\in X$, $\vert x\vert =x\vee (-x)$. An element $x\in E$ is called positive if $x\ge 0$. The positive cone in $E$ is denoted by $E_+$. A Banach lattice $E$ is called an $AL$-space if its norm is additive on the positive cone, i.e., $\Vert f+g\Vert=\Vert f\Vert +\Vert g\Vert$ for all $f,g\in E_+$. An operator $ P\in \calL(E,F)$ is called positive if it maps the positive cone $E_+$ into $F_+$, i.e., $ P E_+\subset F_+$. The set of all such positive operators is denoted by $ \calL_+(E,F) $. An operator $ P\in \mathcal{L}(E,F)$ is called regular if it can be expressed as $ P =P_1- P_2$ for some $ P_{1},P_{2}\in \mathcal{L}_+(E,F) $. The space of all regular operators is denoted by $\mathcal{L}^r(E,F)$.  For further details on Banach lattice theory; see, e.g., \cite{CHARALAMBOS} or \cite{Schaf}.

\section{Admissibility of positive control operator}\footnote{This section refines and completes some results initially presented in \cite[Sect. 2]{El1}.}\label{Sec:2}
In this section, we study the well-posedness and positivity properties of solutions to the abstract differential equation \eqref{diff-equ}. Let $X$ be a Banach lattice. We say that the semigroup $\T$ is positive if $T(t)x\ge 0$ for any $x\in X_+$ and any $t\ge 0$. We say that $A$ is resolvent positive if there is $\omega\in \mathbb{R}$ such that $(\omega,\infty)\subset \rho (A) $ and $R(\la,A)\geq 0$ for all $\la > \omega$. The semigroup $\T$ is called positive if $T(t)x\ge 0$ for any $x\in E_+$ and any $t\ge 0$. The operator $A$ is called resolvent positive if there exists $\omega\in \mathbb{R}$ such that $(\omega,\infty)\subset \rho (A) $ and $R(\la,A)\geq 0$ for all $\la > \omega$. While every generator of a positive C$_0$-semigroup is resolvent positive, the converse is not true in general. Counterexamples and further discussion can be found in, e.g.,  \cite[Sect. 3]{Arendt}. For an overview of the theory of positive C$_0$-semigroups, we refer for instance to the monographs \cite{BFR,Nagel}.
\vspace{.1cm}

We note that $X_{-1,A}$ is not in general a Banach lattice (cf. \cite[Rem. 2.5]{BJVW}), and positivity in $X_{-1,A}$ is defined by the closure of $X_+$ in $X_{-1,A}$. Moreover, we have $X_+=X\cap (X_{-1,A})_+$ (cf. \cite[Prop. 2.3]{BJVW}), where $(X_{-1,A})_+$ denotes the closure of $X_+$ with respect to $\Vert \cdot\Vert_{-1}$, see \cite[Sect. 2.2]{AGPF} for more details. 
\vspace{.1cm}

For $\t\ge 0$, the truncation operator $\mathcal{P}_\t$ and the left shift operator $\calS_\t$ are defined, respectively, by
\begin{align*}
	&(\mathcal{P}_\t u)(t)=
	\begin{cases}
		u(t),&  {\rm for }\;\; t\in [0,\t),\cr
		0,& {\rm for }\;\; t\ge \t,
	\end{cases}\qquad
	(\mathcal{S}_\t u)(t)=
	\begin{cases}
		0,&  {\rm for }\;\; t\in [0,\t),\cr
		u(t+\t),& {\rm for }\;\; t\ge \t.
	\end{cases}
\end{align*}

The following result extends \cite[Prop. 4.3]{Tilman}, establishing a characterization for the well-posedness and positivity of linear infinite-dimensional systems with unbounded control operators.
\begin{proposition}\label{S2.P1} 
	Let $ X$ and $U$ be Banach lattices, $A $ be the generator of a C$_0$-semigroup $ \mathbb{T}$ on $ X $, and $ B\in \mathcal{L}(U,X_{-1,A}) $. Then, the following assertions are equivalent:
	\begin{itemize}
		\item[$(\emph{i})$] For every $x\in X_+$ and $ u\in L^{p}_{loc,+}(\R_+;U) $ the mild solution of \eqref{diff-equ}, given by \eqref{S2.2}, remains in $X_+ $ for all $t\ge 0$.
		\item[$(\emph{ii})$] $\mathbb{T},B $ are positive and, for some $ \tau>0,\,  \Phi_\tau u \in X_+ $ for all $ u\in L^{p}_+([0,\tau];U) $. 
	\end{itemize}
	Moreover, if one of these conditions is satisfied then the differential equation \eqref{diff-equ} admits a unique positive mild solution $ z(\cdot)\in C(\R_+;X)$ given by $$z(t,x,u)=T(t)x+\Phi_tu,$$ for all $t\ge 0$ and $(x,u)\in X\times L^p_{loc}(\R_+,U)$ with  $1\le p<\infty$.
\end{proposition} 	
\begin{proof}
	The implication $ (\emph{ii})\implies (\emph{i}) $ follows directly from \eqref{S2.2}. To prove the converse, let $u$ be identically zero in  \eqref{S2.2}. Then,
	\begin{align*}
		T(t)x=z(t,x,0)\in X_+,
	\end{align*}
	for all $x\in X_+$ and $t\ge 0$, proving that $\T$ is positive. On the other hand, letting $x=0$ in \eqref{S2.2} for arbitrary vector $u\in U_+$ and $t\geq 0$, we obtain
	\begin{align*}
		\tfrac{1}{t}\Phi_t u= \tfrac{1}{t}z(t,0,u) \in X_+,
	\end{align*}
	which tends to $Bu$ as $t\to  0$. The closedness of $X_{-1,+}$ further yields that $Bu\in X_{-1,+}$. Thus, $ (\emph{i})  $ implies $ (\emph{ii}) $.
	\vspace{.1cm}
	
	Now, from $ (\emph{ii}) $ we have $ \Phi_\t u \in X_+ $ for all $u\in L^{p}_+(\R_+;U)$. Since $ \Phi_\t $ is linear and $X$ is a Banach lattice, it follows from \cite[Thm. 1.10]{CHARALAMBOS} that $ \Phi_\t  $ admits a unique extension to a positive operator from $ L^{p}(\R_+;U) $ to $ X $. By Definition \ref{S2.D1}, this implies that $ (A,B) $ is $L^p$-admissible, and thus $\Phi_\tau\in \calL(L^{p}(\R_+;U),X)$ for any $\tau\ge 0$. Therefore, by \cite[Thm. 4.3.1-(i)]{Staf}, the abstract differential equation \eqref{diff-equ} has a unique mild solution $z(\cdot)\in C(\R_+;X)$. This completes the proof.
\end{proof}

We note that $ u $ and $ z $ from \eqref{S2.2} have Laplace transforms related by	
\begin{eqnarray*}
	\hat{z}(\la)=R(\lambda,A)x+\widehat{(\Phi_\cdot u)}(\lambda),\quad \widehat{(\Phi_\cdot u)}(\lambda)=R(\lambda,A_{-1})B \hat{u}(\lambda),
\end{eqnarray*}
for all ${\rm Re}\, \lambda\ge \alpha$, where $\hat{u}$ denote the Laplace transform of $u$ and $\alpha>\omega_{0}(A)$.
\vspace{.1cm}

The following result is a simple consequence of \cite[Prop. 4.3]{AGPF}, which gives characterizations of the positivity of $B$ in terms of the input-maps of $(A,B)$.
\begin{lemma}\label{S2.LB}
	Let $ X$ and $U$ be Banach lattices, and let $ \mathbb{T} $ be a positive $ C_0 $-semigroup on $ X $. Then, the following assertions are equivalent:
	\begin{itemize}
		\item[$(\emph{i})$] $B$ is positive.
		\item[$(\emph{ii})$] $R(\lambda,A_{-1})B$ is positive for all sufficiently large $\lambda\in \R$.
		\item[$(\emph{iii})$] $\Phi_t$ is positive for all $t\geq 0$.
	\end{itemize}
\end{lemma}	
\begin{proof}
	$(\emph{i}) \iff (\emph{ii})$ see \cite[Prop. 4.3]{AGPF}. On the other hand, using the integral formula \eqref{S2.2'} and the fact that
	$$
	B v=\underset{t\to 0 }{\lim}\frac{1}{t}\Phi_t v\qquad ({\rm in}\; X_{-1,A})
	$$
	one readily obtains $(\emph{i})$ $\iff$ $(\emph{iii})$, so that the proof is complete.
\end{proof}

Using linear control systems terminology (see, e.g., \cite[Def. 2.1]{WC}), one can prove that abstract positive control systems are fully determined by their behavior with respect to positive controls. Notably, we have the following result.
\begin{proposition}\label{Admissible-control}
	Let $X$ and $U$ be Banach lattices, $\mathbb{T}$ be a positive $ C_0 $-semigroup, and $p\ge 1$ be a real number. Let $ \Phi:=(\Phi_\t)_{\t\geq 0} $ be a family of maps $ \Phi_\t: L^{p}_+(\R_+;U)\to X_+ $ such that
	\begin{subequations} 
		\begin{align}
			&\Phi_\t (u+v)= \Phi_\t u+ \Phi_\t v\label{eq:1}\\
			&\Phi_{\t+t}u=T(\t)\Phi_t \mathcal{P}_t u+\Phi_\t\mathcal{S}_t u,\label{eq:2}
		\end{align}
	\end{subequations}
	for any $ u,v\in L^{p}_+(\R_+;U)$ and any $ \t,t\geq 0 $. Then, $(T,\Phi)$ is an abstract positive linear control system on $L^{p}(\R_+;U)$ and $X$. In addition, there exists a unique positive control operator $ B\in \mathcal{L}(U,X_{-1,A})$ such that 
	\begin{eqnarray*}
		\Phi_tu =\int_{0}^{t}T_{-1}(t-s)Bu_+(s)ds- \int_{0}^{t}T_{-1}(t-s)Bu_-(s)ds,
	\end{eqnarray*}
	for all $u\in L^{p}(\R_+;U)$ and $t\geq 0 $.
\end{proposition}
\begin{proof}
	It follows from \cite[Thm. 1.10]{CHARALAMBOS} that the map $ \Phi_\t $ has a unique extension to a positive linear operator from $ L^{p}(\R_+;U) $ into $ X $ for each $\t\ge 0$. Further, the extensions, denoted by $\Phi_\t $ again, are given by
	\begin{align}\label{composition}
		\Phi_{\t}u= \Phi_{\t}u_+-\Phi_{\t}u_-,
	\end{align}
	for all $ u\in L^{p}(\R_+;U)$ and $\t\ge 0$. Thus, for every $\t\ge 0$, $ \Phi_\t\in \mathcal{L}(L^{p}(\R_+;U),X) $ as an everywhere defined positive operator is bounded, cf. \cite[Thm. II.5.3]{Schaf}. In addition, it follows from \eqref{eq:2} together with \eqref{composition} that
	\begin{eqnarray*}
		\Phi_{\t+t}u&=&\Phi_{\t+t}u_+-\Phi_{\t+t}u_-\\
		&=&T(\t)\Phi_t \mathcal{P}_t u_++\Phi_\t\mathcal{S}_t u_+ -T(\t)\Phi_t \mathcal{P}_t u_-+\Phi_\t\mathcal{S}_t u_-\\
		&=&T(\t)\Phi_t \mathcal{P}_t(u_+-u_-)+ \Phi_\t\mathcal{S}_t (u_+-u_-)\\
		&=&T(\t)\Phi_t \mathcal{P}_t u+\Phi_\t\mathcal{S}_t u,
	\end{eqnarray*}
	for any $ u,v\in L^{p}(\R_+;U)$ and any $ \t,t\geq 0 $, where we used the linearity of the operators $\mathcal{P}_t$ and $\mathcal{S}_t$ to obtain the last equation. Therefore, in view of \cite[Def. 2.1]{WC}, $(T,\Phi)$ is an abstract linear control system on $L^{p}(\R_+;U)$ and $X$. Furthermore, according to \cite[Thm. 3.9]{WC}, there exists a unique control operator $ B\in \mathcal{L}(U,X_{-1,A})$ such that
	\begin{align*}
		\Phi_tu =\int_{0}^{t}T_{-1}(t-s)Bu(s)ds,
	\end{align*}
	for all $u\in L^{p}(\R_+;U)$ and $t\geq 0 $. Lemma \ref{S2.LB} further yields that $B$ is positive. Additionally, using \eqref{composition}, one obtains
	\begin{align*}
		Bv&= \lim_{\t \mapsto 0}\tfrac{1}{\t}\Phi_{\t}v\\
		&= \lim_{\t \mapsto 0}\tfrac{1}{\t}\Phi_{\t}v_+-\lim_{\t \mapsto 0}\tfrac{1}{\t}\Phi_{\t}v_-\\
		&=Bv_+-Bv_-.
	\end{align*}
	for all $v\in U$. This concludes the proof.
\end{proof}
\begin{remark}
	Note that the above result asserts that positive linear input-maps of abstract positive linear control systems are extensions of families of positive operators that satisfy \eqref{eq:1}-\eqref{eq:2}, but not necessarily the homogeneity property.
\end{remark}
\begin{remark}\label{S2.R2}
	We note that the definition of $L^p$-admissibility of positive control operators for positive semigroups can alternatively be formulated as follows
	\begin{align*}
		\Phi_\t u\in  X_+,\qquad  \forall\, u\in L^p_+(\R_+;U),
	\end{align*}
	for some $\tau > 0$. In this case, we say that $(A, B)$ is positive $L^p$-admissible.
	
	This formulation is clearly equivalent to the standard definition, since every input $u \in L^p(\mathbb{R}_+;U)$ can be written as the difference of two positive inputs. The advantage of this formulation is that admissibility can be deduced from the behavior of the input maps on the positive cone. In particular, in view of \cite[Thm. 1.10]{CHARALAMBOS}, this remains valid even when these maps are defined only on the positive cone and are additive but not necessarily homogeneous.
\end{remark}
Let us now denote by $ \mathfrak{B}_p(U,X,A) $ the vector space of all $L^p$-admissible control operators $B$ for $A$, which is a Banach space for $p\in [1,\infty]$ with the norm 
\begin{align*}
	\Vert B\Vert_{\mathfrak{B}_p}:=\sup_{\Vert u\Vert_{_{L^{p}([0,\t];U)}}\leq 1}\left\Vert \int_{0}^{\tau}T_{-1}(\tau-s)Bu(s)ds\right\Vert,
\end{align*} 
where $\t> 0$ is fixed, cf. \cite[Def. 4.5]{WC}.

\begin{remark}
	Note that $ \mathfrak{B}_{p,+}(U,X,A) $, the set of all $L^p$-admissible positive control operators, forms a positive convex cone in $ \mathfrak{B}_{p}(U,X,A) $. This induces a partial order on $ \mathfrak{B}_{p}(U,X,A) $: for $ B,B'\in  \mathfrak{B}_p(U,X,A)$, we say $ B'\leq B$ if $ B-B'\in \mathfrak{B}_{p,+}(U,X,A) $.
\end{remark}

Next, we establish a sufficient condition for $L^1$-admissibility of control operators for positive semigroups. For this purpose, we recall that $A$ is resolvent positive if there is $\omega\in \mathbb{R}$ such that $(\omega,\infty)\subset \rho (A) $ and $R(\la,A)\geq 0$ for all $\la > \omega$.
\vspace{.1cm}

\begin{theorem}\label{S2.T1}
	Let $ X , U$ be Banach lattices. Assume that $ A$ is a densely defined resolvent positive operator such that, for some $ \la_0>s(A) $ and $c >0$,
	\begin{align}\label{inverse}
		\Vert R(\la_0,A)x\Vert\geq c \Vert x \Vert, \qquad \forall x\in X_+.
	\end{align}
	If $B\in \calL_+(U,X_{-1,A})$, then $(A,B)$ is positive $L^1$-admissible. 
\end{theorem}
\begin{proof}
	We first note that the resolvent inverse estimate \eqref{inverse} implies that $A$ generates a positive $ C_0 $-semigroup $\mathbb{T}$ on $ X $ such that $\omega_{0}(A)=s(A)$, cf. \cite[Thm. 2.5]{Arendt} (see also \cite{BaRo}).
	\vspace{.1cm}
	
	Now, let $ B\in \calL_+(U,X_{-1,A}) $ and $ v  $ be a $U$-valued positive step function, i.e., $$ v(t)=\sum_{k=1}^{n} v_k \1_{[\t_{k-1},\t_{k}[}(t), $$ where $ v_1,\ldots,v_n\in U_+ $, and $ \1_{[\t_{k-1},\t_{k}[} $ denotes the indicator function of $ [\t_{k-1},\t_{k}[$ for $ 0=\t_0<\ldots< \t_n=\t $. Using \cite[Lem. II.1.3-(iv)]{NEN} and the fact that $\T$ and $B$ are positive, one can easily shows that $ \Phi_{\tau}v \in X\cap X_{-1,+}$ and thus $ \Phi_{\tau}v \in X_+$ (cf. \cite[Prop. 2.3]{BJVW}). On the other hand, let $ M\geq 1 $ and $ w> \omega_0(A)$ such that $ \Vert T(t)\Vert \leq Me^{wt} $ for all $ t\geq 0 $. Then, using the resolvent inverse estimate \eqref{inverse}, one obtains that
	\begin{align*}
		c\Vert \Phi_{\tau}v\Vert& \le \left\Vert R(\la_0,A)  \Phi_{\tau}v\right\Vert\\
		&= \left\Vert  \sum_{k=1}^{n}\int_{\t_{k-1}}^{\t_k}T(\tau -s)R(\la_0,A_{-1})B v_kds\right\Vert\\
		&\leq \int_{0}^{\t}\Vert T(\tau -s)R(\la_0,A_{-1})Bv(s)\Vert ds\\
		&\leq  M\Vert R(\la_0,A_{-1})B\Vert  e^{\vert w\vert \t} \Vert v\Vert_{ _{L^{1}([0,\t];U)}},
	\end{align*}
	for any step function $v:[0,\tau]\to U_+$, where we used the fact that $ R(\la_0,A_{-1})B \in \mathcal{L}(U,X) $ (since an everywhere defined positive operator on Banach lattices is norm bounded). Thus,
	\begin{align}\label{Estimate}
		\Vert \Phi_{\tau}v\Vert \le \frac{1}{c}M\Vert R(\la_0,A_{-1})B\Vert  e^{\vert w\vert \t} \Vert v\Vert_{ _{L^{1}([0,\t];U)}},
	\end{align}
	for all positive step function $v$. Therefore, using the fact that $U$-valued positive step functions are dense in $ L^{1}_+([0,\t];U) $ and according to Remark \ref{S2.R2}, it follows from \eqref{Estimate} that $ (A,B) $ is positive $ L^1 $-admissible. This completes the proof.
\end{proof}
\begin{remark}
	We note that Theorem \ref{S2.T1} provide a sufficient conditions for $L^1$-admissibility of positive input operators for positive semigroups, namely a lower norm estimate for the resolvent on the positive cone. Note also that the result of the above theorem is more interesting when the state-space $X$ is not reflexive. In fact, if $X$ is a reflexive Banach lattice, then $B\in \calL(U,X)$, cf. \cite[Thm. 4.8]{WC}.
\end{remark}
\begin{remark}
	It is worth noting that if $A$ has a compact resolvent, then the resolvent inverse estimate \eqref{inverse} implies that the state-space is an $AL$-space, cf. \cite[Cor. B.3]{AGPF}.
\end{remark}


\begin{corollary}\label{CBR}
	Let the assumption of Theorem \ref{S2.T1} be satisfied. If $B\in \calL^r(U,X_{-1,A})$, then $(A,B)$ is $L^1$-admissible.
\end{corollary}
\begin{proof}
	The proof follows from Theorem \ref{S2.T1} and the fact that every regular operator $B$ can be written as the difference of two positive operators.
\end{proof}

For left-invertible positive semigroups, we have the following result. 
\begin{corollary}\label{Cor-left}
	Let $X$ be an $AL$-space and $U$ be a Banach lattice. Assume that $\T$ is positive and left-invertible. Let $B\in \calL^r(U,X_{-1,A})$. Then, $(A,B)$ is $L^1$-admissible.
\end{corollary}
\begin{proof}
	According to \cite[Thm. 2.1]{Xu}, the left-invertibilty of $\T$ implies that there exist $N>0$ and $\alpha>0$ such that
	\begin{align}\label{leftinver}
		\Vert T(t)x\Vert \ge Ne^{-\alpha t}\Vert x\Vert, \qquad \forall x\in X,\, t>0.
	\end{align}
	Using \eqref{leftinver}, the positivity of $\T$, and the additivity of the norm on the positive cone $X_+$, one obtains 
	\begin{align*}
		\Vert R(\la,A)x\Vert&=\left\Vert\int_{0}^{+\infty} e^{-\la t} T(t)x dt\right\Vert\\
		&=\int_{0}^{+\infty} e^{-\la t} \Vert T(t)x\Vert dt\\
		&\ge N \int_{0}^{+\infty} e^{-(\la+\alpha) t} \Vert x\Vert dt=\frac{N}{\la+\alpha}  \Vert x\Vert,
	\end{align*}
	for all $x\in X_+$ and $\la>-\alpha$. Thus, $(A,B)$ is $L^1$-admissible (cf. Corollary \ref{CBR}).
\end{proof}
\begin{remark}
	The above result can be regarded as the counterpart of \cite[Thm. 5.2.2]{TW} in the Banach lattice setting. More precisely, it shows that if the state-space $X$ is an $AL$-space and $\T$ is positive and left-invertible, then every positive control operator $B$ is $L^1$-admissible.
\end{remark}

\begin{example}
	Let us consider the following transport equation
	\begin{align*}
		\begin{cases}
			\frac{\partial}{\partial t}z(t,x)= -\frac{\partial}{\partial x} z(t,x),& t>0, x\in (0,L),\\
			z(t,0)=az(t,L),& t>0,
		\end{cases}
	\end{align*}
	where $a\ge 1$ and $L>0$. On $X:=L^1(0,L)$, we consider the operator $A:D(A)\to X$ defined by 
	\begin{align*}
		Af:=-\frac{\partial}{\partial x}f,\quad
		f\in D(A)=\left\{f\in W^{1,1}(0,L):\; f(0)=af(L)\right\}.
	\end{align*}  
	It is not difficult to see that $A$ is a densely defined resolvent positive operator. We shall show that $A$ satisfies the resolvent inverse estimate \eqref{inverse}. In fact, let $\la >\frac{{\rm ln}(a)}{L}$ and $0\le f\in D(A)$. Set $g=R(\la,A)f$ the positive solution of
	\begin{align*}
		\la f(x)-Af(x)=g(x), \qquad x\in [0,L].
	\end{align*}
	Integrating with respect to $x$, we obtain
	\begin{align*}
		\la \Vert f\Vert_X+ f(L)-af(L)=\Vert g\Vert_X.
	\end{align*}
	Since $a\ge 1$, then
	\begin{align*}
		\Vert R(\la,A)f\Vert_X\ge \frac{1}{\la}\Vert f\Vert_X, \qquad \forall f\in X_+.
	\end{align*}
\end{example}
\begin{example}
	As discussed in \cite[Appendix B]{AGPF}, a notable class of semigroups that satisfy the resolvent inverse estimate \eqref{inverse} consists of Markov semigroups in $L^1$-spaces. In fact, let $(\Omega,\Sigma, m)$ be a $\sigma$-finite measure spaces and denote by $L^1:=L^1(\Omega,\Sigma, m)$ the corresponding space of integrable functions. We say that $ \M:=(M(t))_{t\ge 0}$ is a Markov semigroup in $L^1$-space if, for every $t\ge 0$ and $f\in L_+^1$, 
	\begin{itemize}
		\item[(\emph{i})] $M(t)f\ge 0$;
		\item[(\emph{ii})] $\Vert M(t)f\Vert_{L^1}=\Vert f\Vert_{L^1}$.
	\end{itemize}
	\vspace{.1cm}
	
	Let $\M$ a Markov semigroup in $L^1$ and denote by $Q$ its corresponding generator. Then, using the additivity of $L^1$-norm on the positive cone, one obtains  
	\begin{align*}
		\Vert R(\la,Q)f\Vert_{L^1}=\int_{0}^{+\infty}e^{-\la t}\Vert M(t)f\Vert_{L^1}dt=\frac{1}{\la} \Vert f\Vert_{L^1},
	\end{align*}
	for all $f\in X_+$ and $\la>s(Q)\vee 0$. Thus, for $\la_1\ge \la_0>s(Q)\vee 0$, one has 
	\begin{align*}
		\Vert R(\la_0,Q)f\Vert_{L^1}\ge \Vert R(\la_1,Q)f\Vert_{L^1}=\frac{1}{\la_1} \Vert f\Vert_{L^1}, \qquad \forall f\in L_+^1.
	\end{align*}
	For concrete examples of such semigroups, particularly those arising in applications to biological systems and stochastic processes, we refer to \cite{Rudnicki}.
\end{example}

\section{Robustness of ISS}\label{Sec:4}
In this section, we examine the invariance of input-to-state stability under Desch-Schappacher perturbations. More precisely, we consider the perturbed linear control system
\begin{align}\label{Perturbed}
	\begin{cases}
		\dot{z}(t)=(A_{-1}+P)\vert_{X}z(t)+Bu(t), & t>0,\\
		z(0)=x,&
	\end{cases}
\end{align}
where $A$ generates a C$_0$-semigroup $\T$ on $X$, $P\in \calL(X,X_{-1,A})$, and $B\in \calL(U,X_{-1,A})$. Here, $(A_{-1}+P)\vert_{X}:D((A_{-1}+P)\vert_{X})\to X_{-1,A}$ is the operator defined by 
\begin{align*}
	&(A_{-1}+P)\vert_{X}:= A_{-1}+P,\\
	& D((A_{-1}+P)\vert_{X})=\left\{f\in X:\; (A_{-1}+P)f\in X\right\}.
\end{align*}
To proceed, we introduce key concepts. A semigroup $\T$ is strongly exponential stable, if there exist $M,\theta> 0$ such that
\begin{align*}
	\Vert T(t)x\Vert_{X} \leq Me^{-\theta t}\Vert x\Vert_{D(A)}, \qquad \forall t\ge 0,\, x\in D(A).
\end{align*}
If the right-hand side of the above inequality holds with respect to the norm $\Vert\cdot\Vert_X$, then $\T$ is uniform exponential stable, which is equivalent to $\omega_0(A)<0$. Since $s(A)\le \omega_0(A)$, uniform exponential stability implies $s(A)<0$. When $\omega_0(A)=s(A)$, one says that $\mathbb{T} $ satisfies the spectrum determined growth property (\cite[p. 161]{CZ}). This property holds for many classes of semigroups, such as positive semigroups in $L^p$-space for $p\in [1,+\infty)$; see, e.g., \cite{Lutz1} or \cite[Thm. 1]{Lutz2}.
\vspace{.2cm}

We now recall the definition of ISS for linear control systems. Indeed, according to \cite[Def. 2.7]{JaNPS} (see also \cite[Def. 3.17]{Andri}), we select the following definition.
\begin{definition}\label{ISS-def}
	Let $A$ be the generator of a C$_0$-semigroup $\T$ on $X$ and $p\ge 1$ a real. Then, system \eqref{diff-equ} is called exponentially $L^p$-input-to-state ($eL^p$-ISS for short) with a linear gain if there exist $N,\mu,G>0$ such that for every $t\ge 0$, $x\in X$ and $u\in L^p(\R_+;U)$
	\begin{itemize}
		\item[(\emph{i})] $z(t)\in C(\R_+;X)$,
		\item[(\emph{ii})] $
		\Vert z(t,x,u)\Vert_X\le Ne^{-\mu t}\Vert x\Vert_X+G\Vert u\Vert_{L^p(\R_+;U)}.
		$
	\end{itemize}
\end{definition}
\begin{remark}
	It follows from the above definition that e$L^p$-ISS of \eqref{diff-equ} is fully characterized by the uniform exponential stability of $\T$ and the $L^p$-admissibility of $(A,B)$ for $p\in [1,+\infty)$. In this case, $G:=\sup_{t>0}\kappa(t)$ where $\kappa(\cdot)$ is the constant of admissibility in \eqref{S2.2'}, see \cite[Rem. 2.12]{JaNPS}. If, additionally, $z(\cdot)$ is continuous w.r.t. time in the norm of $X$, then \eqref{diff-equ} is exponentially $L^\infty$-ISS (see, \cite[Sect. 2]{JaNPS} and \cite[Sect. 3.2]{Andri} for more details).
\end{remark} 
\begin{remark}
	We note that the estimate in Definition \ref{ISS-def} is a special variant of ISS estimates. Such stability estimates have proven to be a powerful tool to simultaneously investigate the internal stability and the robustness of dynamical systems to external controls or disturbances. For an overview of recent developments in ISS, the reader is referred to the survey paper \cite{Andri} and the books \cite{KaKr,Mir}.
\end{remark}

We state the following proposition.
\begin{proposition}\label{Generation-prop}
	Let $A$ generate a positive C$_0$-semigroup $\T$ on a Banach lattice $X$ and $P\in \calL_+(X,X_{-1,A})$. Assume that $(A,P)$ is positive $L^p$-admissible for some $p\in [1,+\infty)$. Then,
	\begin{itemize}
		\item[(\emph{i})] $(A_{-1}+P)\vert_X$ generates a positive $ C_0 $-semigroup $\mathbb{S}:=(S(t))_{t\ge 0}$ on $ X $ with $0\le \T\le \S$.
		\item[(\emph{ii})] We have $s(A)\le s((A_{-1}+P)\vert_X)$ and $R(\la,A)\le R(\la, (A_{-1}+P)\vert_X)$ for all $\la >s((A_{-1}+P)\vert_X) $.
	\end{itemize}
\end{proposition}
\begin{proof}
	(\emph{i})- For $\t_0>0$ and $p\in [1,+\infty)$, we denote by $\Psi_{\t_0}\in \calL(L^p([0,\tau_0];X),X_{-1,A})$ the operator defined by
	\begin{align*}
		\Psi_{\t_0} f=\int_{0}^{\tau_0}T_{-1}(\tau_0-s)Pf(s)ds,
	\end{align*}
	for $f\in L^p([0,\tau_0];X)$. Since $(A,P)$ is positive $L^p$-admissible, then $\Psi_{\t_0}f\in X_+$ for all $f\in  L^p_+([0,\tau_0];X)$. Using the linearity of $\Psi_{\t_0}$ and the fact that $X_+,L^p_+([0,\tau_0];X)$ are generating cones, one has ${\rm Range}(\Psi_{\t_0}) \subset X$. Thus, according to \cite[Cor. II.3.4]{Nagel}, the operator $(A_{-1}+P)\vert_X$ generates a $ C_0 $-semigroup $\mathbb{S}:=(S(t))_{t\ge 0}$ on $ X $ given by 
	\begin{align}\label{variation1}
		S(t)x=T(t)x+\int_{0}^{t}T_{-1}(t-s)PS(s)xds,
	\end{align}
	for all $t\ge 0$ and $x\in X$. The Dyson-Phillips expansion formula (see, e.g., \cite[Cor. III.3.2]{Nagel}) further yields that $0\le T(t)\le S(t)$ for all $t\ge 0$. This shows (\emph{i}).
	
	(\emph{ii})- Using (\emph{i}) and \cite[Thm. 12.7]{BFR}, one obtains 
	\begin{align*}
		0\le R(\la,A)\le R(\la, (A_{-1}+P)\vert_X),
	\end{align*}
	for all $\la >\max\{s(A),s( (A_{-1}+P)\vert_X)\}$. Thus, for such $\la$, we have 
	\begin{align*}
		\Vert R(\la,A)\Vert\le \Vert R(\la, (A_{-1}+P)\vert_X)\Vert.
	\end{align*}
	Therefore, we obtain $s(A)\le s((A_{-1}+P)\vert_X)$, since $s(A)\in \sigma(A)$ (see, e.g., \cite[Cor. 12.9]{BFR}). This proves (\emph{ii}). The proof is complete.
\end{proof}

Using the above proposition and according to \cite[Thm. 3.2]{BJVW}, one can deduce the following stability result.
\begin{lemma}\label{Exp-stab}
	Let $A$ be the generator of a positive C$_0$-semigroup $\T$ on a Banach lattice $X$ and $P\in \calL_+(X,X_{-1,A})$. Then, 
	\begin{align*}
		s((A_{-1}+P)\vert_X)<0 \iff s(A)<0\; {\rm and }\; r(R(0,A_{-1})P)<1.
	\end{align*} 
\end{lemma}
\begin{proof}
	The proof follows from \cite[Thm. 3.2]{BJVW}. For the sake of completeness we include here the details.
	
	If $s((A_{-1}+P)\vert_X)<0$, then, by Proposition \ref{Generation-prop}-(ii), we obtain $s(A)<0$. Moreover, using the positivity of $\S$, we obtain $0\in \rho((A_{-1}+P)\vert_X)$ and $R(0, ((A_{-1}+P)\vert_X)) \ge 0$. Thus, according to \cite[Thm. 3.2]{BJVW}, we have $r(R(0,A_{-1})P)<1$. 
	
	Conversely, assume that $s(A)<0$ and $r(R(0,A_{-1})P)<1$. Let $s((A_{-1}+P)\vert_X)>-\infty$, as otherwise there is noting to prove. It follows from \cite[Thm. 3.2]{BJVW} that $0\in \rho((A_{-1}+P)\vert_X)$ and $R(0, ((A_{-1}+P)\vert_X)) \ge 0$, so that $s((A_{-1}+P)\vert_X)<0$. This completes the proof.
\end{proof}

As a consequence, we obtain the following spectral characterization for the exponential input-to-state stability of the perturbed system \eqref{Perturbed}.
\begin{proposition}\label{Main-charac}
	Let the assumptions of Proposition \ref{Generation-prop} be satisfied and $B\in \calL(U,X_{-1,A})$. Assume that $\S$ satisfies the spectrum determined growth property and $(A,B)$ is $L^q$-admissible for some $q\in [1,\infty)$. Then, for any $r\ge q$, system \eqref{Perturbed} is $eL^r$-ISS if and only if 
	\begin{align*}
		s(A)<0 \quad {\rm and }\quad r(R(0,A_{-1})P)<1.
	\end{align*}
\end{proposition}
\begin{proof}
	Since $(A,P)$ is positive $L^p$-admissible for some $p\in [1,\infty)$, we know from Proposition \ref{Generation-prop} that $(A_{-1}+P)\vert_X$ generates a positive $ C_0 $-semigroup $\S$ on $X$. On the other hand, according to \cite[Thm. 3.3-(ii)]{Hadd} (see also \cite[Cor. 2.7]{Mei}), we have 
	\begin{align}\label{vector-spqce-equality}
		\mathfrak{B}_{q}(U,X,A) = \mathfrak{B}_{q}(U,X,(A_{-1}+P)\vert_{X}), \quad \forall q\in [1,\infty).
	\end{align} 
	If $(A,B)$ is $L^q$-admissible for some $q\in [1,\infty)$, then $((A_{-1}+P)\vert_{X},B)$ is also $L^q$-admissible. Thus, by H\"{o}lder’s inequality $((A_{-1}+P)\vert_{X},B)$ is $L^r$-admissible for any $r\ge q$ and zero-class admissible for any $r>q$ (i.e., $\Vert \Phi_\t^{(A_{-1}+P)\vert_{X}}\Vert\to 0$ as $\t\to 0$). Therefore, using Lemma \ref{Exp-stab} and the fact that $\S$ satisfies the spectrum determined growth property, it follows from \cite[Prop. 2.10]{JaNPS} that system \eqref{Perturbed} is $eL^r$-ISS if and only if $s(A)<0$ and $r(R(0,A_{-1})P)<1$. The proof is completes.
\end{proof}
\begin{remark}
	We note that for $r=\infty$, $((A_{-1}+P)\vert_{X},B)$ is zero-class $L^\infty$-admissible, and by \cite[Prop. 2.5]{JaNPS}, the mild $z(\cdot,x,u)$ is continuous in $t$ for any $f\in X$ and $u\in L^\infty_{loc}(\R_+,U)$. 
\end{remark}

We can now state and prove the main result of this section.
\begin{theorem}\label{F-Main1}
	Let the assumptions of Proposition \ref{Main-charac} be satisfied and $q\in [1,\infty)$. Then, the following assertions are equivalent:
	\begin{itemize}
		\item [(\emph{i})] System \eqref{Perturbed} is $eL^q$-ISS.
		\item[(\emph{ii})] System \eqref{diff-equ} is $eL^q$-ISS and $r(R(0,A_{-1})P)<1$.
	\end{itemize}
\end{theorem}
\begin{proof}
	First, we note that item (\emph{i}) in Proposition \ref{Generation-prop} implies that $\omega_{0}((A_{-1}+P)\vert_X)\ge \omega_{0}(A)$. Using Lemma \ref{Exp-stab} and the fact that $\S$ satisfies the spectrum determined growth property, one easily prove that
	\begin{align}\label{char1}
		\omega_{0}((A_{-1}+P)\vert_X)<0 \iff \omega_{0}(A)<0\; {\rm and }\; r(R(0,A_{-1})P)<1.
	\end{align} 
	Now, using \eqref{vector-spqce-equality} and \eqref{char1}, it follows from \cite[Prop. 2.10]{JaNPS} that
	\begin{align*}
		{\rm System} \, \eqref{Perturbed}\, {\rm is}\, eL^q{\rm -ISS}&\iff 	
		(A_{-1}+P)\vert_X\, {\rm is}\, L^q{\rm -admissible} \, {\rm and}\, \omega_{0}((A_{-1}+P)\vert_{X})<0\\
		&\iff (A,B)\, {\rm is}\, L^q{\rm- admissible},\,\omega_{0}(A)<0,\, {\rm and}\, r(R(0,A_{-1})P)<1\\
		&\iff {\rm System}\, \eqref{diff-equ}\,{\rm is}\, eL^q{\rm -ISS}\,{\rm and }\, r(R(0,A_{-1})P)<1,
	\end{align*}
	where we used the fact that $X_{-1,A}\cong X_{-1, (A_{-1}+P)\vert_{X}}$. This completes the proof. 
\end{proof}

\begin{remark}
	The above result shows that the perturbed system \eqref{Perturbed} is $eL^q$-ISS if and only if (i) the unperturbed system \eqref{diff-equ} is $eL^q$-ISS, and (ii) the small-gain condition $r(R(0, A_{-1})P) < 1$ holds. The condition $r(R(0, A_{-1})P) < 1$ ensures that the perturbation does not amplify instability, providing a precise and necessary criterion for robust stability in the presence of disturbances.
\end{remark}

As a matter of fact, we obtain the following characterization for the eISS of a class of linear positive systems.
\begin{proposition}\label{Invariance-theo}
	Let $X,U$ be Banach lattices, $ A$ a densely defined resolvent positive operator, and $p\ge 1$ be a real number. Assume that:
	\begin{itemize}
		\item[(\emph{i})] $A$ satisfies the resolvent inverse estimate \eqref{inverse} for some $\la_0>0$ sufficiently large;
		\item[(\emph{ii})]  $ P\in \mathcal{L}_+(X,X_{-1,A} )$ and $ B\in \mathcal{L}^r(U,X_{-1,A} )$.
	\end{itemize}
	Then,  
	\begin{align*}
		{\rm System} \; \eqref{Perturbed}\; {\rm is}\; eL^p{\rm -ISS}\iff 	
		s(A)<0\;{\rm and }\; r(R(0,A_{-1})P)<1.
	\end{align*}
\end{proposition}
\begin{proof}
	We know from Corollary \ref{CBR} that $(\emph{i})$ implies that $(A,B)$ is $L^p$-admissible and $\omega_0(A)=s(A)$. On the other hand, using $(\emph{i})$ and the fact that $ P\in \mathcal{L}_+(X,X_{-1,A} )$, it follows from Theorem \ref{S2.T1} that there exists $\t>0$ such that
	\begin{align*}
		\int_{0}^{\t} T(\t-s)P f(s) ds\in X_+, \qquad \forall f\in L^p_+([0,\t];X).
	\end{align*}
	Then, according to Proposition \ref{Generation-prop}, the operator $(A_{-1}+P)\vert_X$ generates a positive $ C_0 $-semigroup $\mathbb{S}$ on $ X $. In addition, we have 
	\begin{align*}
		\Vert R(\la_0, (A_{-1}+P)\vert_X)f\Vert\ge \Vert R(\la_0,A)f\Vert\ge c\Vert f\Vert, \qquad \forall f\in X_+,
	\end{align*}
	since $A$ satisfies the resolvent inverse estimate \eqref{inverse} for some $\la_0$ sufficiently large.  Thus, according to \cite[Thm. 2.5]{Arendt}, $\S$ satisfy the spectrum determined growth property. Hence, by virtue of Proposition \ref{Main-charac}, system \eqref{Perturbed} is exponentially $L^p$-ISS if and only if $s(A)<0$ and $r(R(0,A_{-1})P)<1$.  This completes the proof.
\end{proof}

\section{Application}\label{Sec:5}
We consider the following renewal equation
\begin{align}\label{system}
	\begin{cases}
		\dot{z}(t)= -\displaystyle\frac{\partial}{\partial x} z(t,x)-q(x)z(t,x),& t> 0, \;x >0,\cr  
		z(0,x)= f_0(x), & x>0,\cr
		z(t,0)= \displaystyle\int_{0}^{+\infty}\beta(x)z(t,x) dx+ u(t), & t> 0.
	\end{cases}
\end{align}
Here, the function $z(t,x)$ represents the population density of certain species at time $t$ with size $x$. The functions $\beta\in  L^{\infty}_+(\R_+)$ and $q\in  L^{\infty}_{+}(\R_+)$ denote the birth rate and death rate, respectively. The birth process is determined by the formula $\int_{0}^{\infty} \beta zdx$,
which describes the recruitment process in the population of newborn individuals. The function $u$ represents the birth control. It is important to note that in recent years there has been a growing interest in controllability properties under positivity constraints of age- or size-structured population models; see, e.g., \cite{EFS,HMT} and references therein. 
\vspace{.1cm}

The goal of this section is to analyze an exponential ISS result for \eqref{system}. For this purpose, we consider the state-space $L^1(\R_+)$. Note that the function space $ L^1(\R_+)$ is the natural state space for system \eqref{system} since $L^1$-norm gives the total population size.
\vspace{.1cm}

On $L^1(\R_+)$ we consider the operator $A_m:D(A_m)\to L^1(\R_+)$ defined by 
\begin{eqnarray*}
	A_m:= -\frac{\partial}{\partial x}f-q(\cdot)f,\quad f\in D(A_m):=W^{1,1}(\R_+),
\end{eqnarray*}
We select the following operator
\begin{align*}
	\Gamma f:=f(0), \qquad f\in D(A_m) .
\end{align*}
Using the introduced notation, one can rewrite \eqref{system} as follows
\begin{align}\label{Bounadry-system}
	\begin{cases}
		\dot{z}(t)=A_m z(t),& t>0,\\
		\Gamma z(t)=\Lambda z(t),& t>0,
	\end{cases}
\end{align}
where $z(t):=z(t,\cdot)$ and $\Lambda f:= \int_{0}^{+\infty}\beta(x)f(x) dx$ for all $f\in L^1(\R_+)$. It is clear that $\Lambda\in \calL_+(L^1(\R_+),\R)$.
\vspace{.1cm}

Now, let us consider the operator $A:D(A)\to L^1(\R_+)$ defined by 
\begin{eqnarray*}
	A:=A_m,\quad D(A):=\ker (\Gamma)=\left\{f\in D(A_m):\; f(0)=0\right\}.
\end{eqnarray*}
It is well-known that $A$ generates the right-shift semigroup $\T$ on $L^1(\R_+)$ given, for every $t\ge 0$ and $f\in L^1(\R_+)$, by 
\begin{align*}
	(T(t)f)(x):=\begin{cases}
		e^{-\int_{x-t}^{x}q(\sigma)d\sigma}f(x-t),& t\le x\\
		0,& t>x.
	\end{cases}
\end{align*}
On the other hand, it is clear that the operator $\Gamma$ is onto. It follows from \cite[Lem. 1.2]{Gr} that the Dirichlet operator, associated with $(A_m,\Gamma)$,
\begin{align}\label{Dirichlet}
	D_\lambda:=\left(\Gamma\vert_{\ker(\lambda I_X-A_m)}\right)^{-1}, \quad  \lambda\in \C,
\end{align} 
exists and bounded. Moreover, by a simple calculation, we obtain
\begin{align*}
	(D_\la b)(x)=e^{-\int_{0}^{x}q(\sigma)d\sigma-\la x}b,
\end{align*}
for all $\la\in \C$, $x\ge 0$ and $b\in \R$. Thus, one can verify that the operator
\begin{align}\label{boundary-control}
	B:=(\lambda I_X-A_{-1})D_\lambda\in\mathcal{L}(\R,(L^1(\R_+))_{-1,A}),
\end{align}
is independent of the choice of $\lambda$, Range$(B)\cap L^1(\R_+)=\{0\}$, and that 
\begin{align}\label{representation}
	(A_m-A_{-1})|_{D(A_m)}=B \Gamma.
\end{align}
Using \eqref{representation}, system \eqref{system} becomes
\begin{eqnarray*}
	\dot{z}(t)=(A_{-1}+P)z(t)+Bu(t),\quad t>0,\qquad z(0)=f_0,
\end{eqnarray*}
where $P:=B\Lambda$ and $B:= (-A_{-1})D_0$. From Lemma \ref{S2.LB} it follows that $B$ is positive, and therefore $P$ is also positive. 
\vspace{.1cm}

Let us now prove that $A$ satisfies the resolvent inverse estimate \eqref{inverse}. Using the additivity of the $L^1$-norm on the positive cone, one can easily shows that
\begin{align}\label{Estim-semig}
	\Vert T(t)f\Vert_{L^1(\R_+)}\ge  e^{-\Vert q\Vert_{\infty}t}\Vert f\Vert_{L^1(\R_+)},
\end{align}
for all $t\ge 0$ and $f\in L^1_+(\R_+)$. Thus, using \eqref{Estim-semig} and the additivity of the $L^1$-norm on the positive cone, one obtains
\begin{align*}
	\Vert R(\la, A)f\Vert_{L^1(\R_+)}&=\left\Vert \int_{0}^{+\infty}e^{-\la t}T(t)f dt\right\Vert_{L^1(\R_+)}\\
	&=\int_{0}^{+\infty}e^{-\la t}\left\Vert T(t)f\right\Vert_{L^1(\R_+)}dt
	\ge \frac{1}{\la+\Vert q\Vert_{\infty}}\Vert f\Vert_{L^1(\R_+)},
\end{align*}
for all $\la >-\Vert q\Vert_{\infty}$ and $f\in L^1_+(\R_+)$.
Therefore, 
\begin{align}\label{Inverse-A}
	\Vert R(\la, A)f\Vert_{L^1(\R_+)}\ge \frac{1}{\la+\Vert q\Vert_{\infty}}\Vert f\Vert_{L^1(\R_+)},
\end{align}
for all $\la >-\Vert q\Vert_{\infty}$ and $f\in L^1_+(\R_+)$. Using \eqref{Inverse-A} and the fact that $s(A)<0$, it follows from Proposition \ref{Invariance-theo} that system \eqref{system} is exponentially $L^p$-ISS if and only if $r(R(0,A_{-1})P)<1$. Notice that  
\begin{align*}
	R(0,A_{-1})P=R(0,A_{-1})B\Lambda=D_0\Lambda. 
\end{align*}
Then, the spectral radius $r(R(0,A_{-1})P)$ can be estimated by
\begin{align*}
	r(R(0,A_{-1})P)\le \Vert D_0\Lambda \Vert_{ \calL(L^1(\R_+))}\le \frac{1}{\Vert q\Vert_{\infty}} \Vert \beta\Vert_\infty,
\end{align*}
Hence, for $ \Vert \beta\Vert_\infty<\Vert q\Vert_{\infty}$, we obtain $r(R(0,A_{-1})P)<1$, so that system \eqref{system} $eL^p$-ISS.

\section{Conclusion}
In this paper, we have discussed well-posedness, positivity and stability property of linear distributed parameter systems. Our contribution is twofold: we have presented a resolvent condition for $L^1$-admissibility of control operators. Then, we have established necessary and sufficient conditions for the robustness of exponential ISS under positive Desch-Schappacher perturbations. The practical relevance of our results is demonstrated through an application to a boundary-controlled transport equation with non-local boundary conditions. Future work will extend our results to general boundary control systems.

\end{document}